\documentclass[12pt,leqno]{amsart}

\usepackage{amsthm}
\usepackage[usenames]{color}

\usepackage{amsmath}
\usepackage{amscd}
\usepackage{amsmath}
\usepackage{amssymb}
\usepackage{latexsym}
\makeatletter

\@addtoreset{equation}{section}
\makeatother

\newtheorem{thm}{Theorem}[section]
\newtheorem{pr}[thm]{Proposition}
\newtheorem{lm}[thm]{Lemma}
\newtheorem{cor}[thm]{Corollary}

\theoremstyle{definition}

\newtheorem{df}[thm]{Definition}

\newcommand{\sm}{\raisebox{2.33pt}{~\rule{6.4pt}{1.3pt}~}}

\input xy \xyoption{all} \CompileMatrices

\begin{document}

\title{Holonomic \'etale sheaves are constructible} 

\author{Ahmed Abbes and Takeshi Saito
}

\address{Laboratoire Alexander Grothendieck, UMR 9009 du CNRS, Institut
des Hautes \'Etudes Scientifiques, 35 route de Chartres, 91440
Bures-sur-Yvette, France}
\email{abbes@ihes.fr}

\address{School of Mathematical Sciences,
University of Tokyo,
3-8-1 Komaba,
Meguro, Tokyo 153-8914 Japan}
\email{t-saito@ms.u-tokyo.ac.jp}

\maketitle

\begin{abstract}
Building on Beilinson's work, 
``constructible sheaves are holonomic,'' we introduce the notion of holonomicity for \'etale sheaves, without assuming a priori constructibility. 
We establish the converse of Beilinson's result, showing that holonomic sheaves are indeed constructible. This can be seen as an \'etale analogue of Kashiwara's theorem on holonomic 
${\mathcal D}_X$-modules.
\end{abstract}

For a coherent ${\mathcal D}_X$-module 
${\mathcal M}$
on a complex manifold $X$,
its singular support 
$SS{\mathcal M}$
is defined as a closed conical subset
of the cotangent bundle \cite{SKK},
\cite{HTT}.
Inspired by the analogy between
${\mathcal D}_X$-modules
on complex manifolds
and \'etale sheaves
on algebraic varieties,
Beilinson proved that
every constructible sheaf ${\mathcal F}$
on a smooth scheme $X$ over a field
has a singular support $SS{\mathcal F}$
and that its irreducible components
have the same dimension as
$X$ in \cite{SS}.
This means that
constructible sheaves are holonomic
as the title of his article says
although he did not explicitly
define
the holonomicity.

Following his construction,
we introduce the holonomicity
for \'etale sheaves, 
without assuming a priori that
they are constructible.
We prove conversely
that
holonomic sheaves are
in fact constructible.
One may consider this as an \'etale
analogue of a theorem of Kashiwara \cite{Kas}
asserting that for a holonomic
${\mathcal D}_X$-module 
${\mathcal M}$
on a complex manifold $X$,
the ${\mathbb C}_X$-module
${\mathcal E}{\it xt}^i_X({\mathcal M},{\mathcal O}_X)$
is constructible for every $i$,
which is an important ingredient
of the Riemann--Hilbert correspondence.
We expect a similar notion
of holonomicity
for \'etale sheaves on rigid varieties
could be defined
and would give a good substitute
for the deficient notion of constructibility
in this context.

The contents of the sections are summarized as follows. 
In the brief Section 1, 
we recall the fact that a closed conical subset of a vector bundle is
determined by its base and its projectivization. 
In Sections 2 and 3, 
we study two properties of morphisms 
relative to a closed conical subset $C$ 
of the cotangent bundle, 
largely following \cite{SS}.
The notion of $C$-transversality applies when $C$ 
is a subset for the target of a morphism 
$h\colon W \to X$, 
while $C$-acyclicity applies when $C$ is a subset 
for the source of a morphism 
$f\colon X \to Y$. 
Both notions enter the definition of the micro support of 
a sheaf ${\mathcal F}$; 
in \cite{SS}, they are both referred to as transversality. 
In Section 4, we extend the definition of micro support given in 
\cite{SS} to sheaves that are not necessarily constructible; 
in \cite{SS}, constructibility is assumed. 
In Section 5, we prove the existence of the singular support of 
${\mathcal F}$, defined as the smallest closed conical subset on which 
${\mathcal F}$ is micro supported, 
without assuming constructibility. 
The proof follows that of \cite{SS} and relies on the Radon transform.
Finally, in Section 6, we prove the main result: 
if $\dim SS{\mathcal F} \leqq \dim X$ and 
if every fiber of ${\mathcal F}$ is constructible, 
then ${\mathcal F}$ itself is constructible. 
The proof proceeds by induction on 
$\dim X$ and ultimately reduces to the fact that 
a sheaf micro supported on the zero section is locally constant.

The authors thank 
Tong Zhou 
for, after an earlier version of
this article, asking if
the singular support exists
without the constructibility assumption
and the anonymous referee
for suggesting removing the perfectness assumption
in an earlier version.
This work was done during the second
author's stay at IHES in September 2024.
He thanks the first author 
and the institute for the
hospitality.
The work was partially supported
by Kakenhi 24K06683.

\section{Closed conical subsets}

\begin{df}\label{dfconical}
Let $X$ be a scheme 
and $V$ be a vector bundle on $X$.
We say that a closed subset $C$ of $V$
is conical if
it is stable under the action of
the multiplicative group ${\mathbf G}_m$.

The intersection $B$
of a closed conical subset $C$
with the $0$-section,
regarded as a closed
subset of $X$,
is called the base of $C$.
The projectivization
${\mathbf P}(C)
\subset {\mathbf P}(V)$
is defined to be the image of
the complement
$C\sm B$
by the canonical projection
$V\sm X\to {\mathbf P}(V)$
from the complement
of the $0$-section. 
\end{df}

\begin{lm}\label{lmbase}
Let $X$ be a scheme 
and $C$ be a closed conical subset of 
 a vector bundle $V$ on $X$.

{\rm 1.}
$C$ equals the union of the base
$B$ and the inverse image
of the projectivization
${\mathbf P}(C)$
by $V\sm X\to {\mathbf P}(V)$.

{\rm 2.}
$U=\{x\in X\mid C_x=\varnothing\}
\subset
U'=\{x\in X\mid C_x\subset \{0\}\}$
are open subsets of $X$.
\end{lm}

\proof{
1.
Since $C$ is conical,
$C\sm B$ equals 
the inverse image
of ${\mathbf P}(C)$.

2.
$U$ is the complement of
$B$ and $U'$
is the complement of the image
of 
${\mathbf P}(C)$
by the projection
${\mathbf P}(V)\to X$.
\qed

}
\medskip

\section{$C$-transversality}

Let $k$ be a field and $X$ a smooth scheme over $k$.
The covariant vector bundle $T^*X$ associated to
the locally free ${\mathcal O}_X$-module 
$\Omega^1_{X/k}$
is called the cotangent bundle
of $X$. For a smooth subscheme $Z$ of $X$,
the conormal bundle $T^*_ZX$
associated to the conormal sheaf $N_{Z/X}$
is canonically identified with 
${\rm Ker}(T^*X\times_XZ
\to T^*Z)$.
In particular for $Z=X$,
the $0$-section of $T^*X$ is denoted
by $T^*_XX$.

\begin{lm}\label{lmU}
Let $X$ be a smooth scheme
over $k$.
Let $C\subset T^*X$
be a closed conical subset
satisfying $\dim C\leqq
\dim X$.
Then there exists a dense
open scheme $U\subset X$ 
such that the restriction
$C_U$ is a subset of 
the $0$-section $T^*_UU$.

\end{lm}

\proof{
Since $\dim C\leqq \dim X$,
the projectivization
${\mathbf P}(C)$ has dimension
strictly less than $\dim X$.
Hence it follows from Lemma \ref{lmbase}.
\qed

}

\begin{df}\label{dfCtrans}
{\rm (\cite[1.2]{SS})}
Let $h\colon W\to X$
be a morphism
of smooth schemes over a field $k$
and $C$ be a closed conical
subset of the cotangent bundle
$T^*X$.

{\rm 1.}
Let $w$ be a point of $W$.
We say that
$h$ is $C$-transversal at $w$,
if the intersection of
the inverse image
$C_w\subset T^*X\times_Xw$
of $C$
with the kernel
${\rm Ker}(T^*X\times_Xw
\to T^*W\times_Ww)$ inside 
$T^*X\times_Xw$
is a subset $\{0\}$.

{\rm 2.}
We say that $h$ is $C$-transversal
if $h$ is $C$-transversal
at every point of $W$.
\end{df}

The condition that
$h$ is $C$-transversal means that
the intersection of
the inverse image
$h^*C\subset T^*X\times_XW$
of $C$
with the kernel
${\rm Ker}(T^*X\times_XW
\to T^*W)$ inside 
$T^*X\times_XW$
is a subset of
the $0$-section.

\begin{lm}\label{lmCtrans}
Let $h\colon W\to X$
be a morphism
of smooth schemes over a field $k$
and $C$ be a closed conical
subset of the cotangent bundle
$T^*X$.

{\rm 1. (\cite[Example 1.2 (i)]{SS})}
$h$ is smooth if and only if
$h$ is $T^*X$-transversal.

{\rm 2.}
If $C$ is a subset of
the $0$-section,
then $h$ is $C$-transversal.

{\rm 3.}
If $C'$ is a closed conical subset
of $T^*X$ containing $C$ as a subset
and if $h$ is $C'$-transversal,
then $h$ is $C$-transversal.
Consequently,
if $h$ is smooth,
then $h$ is $C$-transversal.

{\rm 4. (\cite[Lemma 1.2 (i)]{SS})}
The subset $\{w\in W\mid
h$ is $C$-transversal at $w\}$ 
is an open subset of $W$.

{\rm 5.}
Let $Z\subset X$ be a closed subscheme
smooth over $k$
and assume that $C$ is
the conormal bundle $T^*_ZX\subset T^*X$.
Then the following conditions are equivalent:

{\rm (1)}
$h$ is $C$-transversal.

{\rm (2)}
$V=Z\times_XW$ is smooth over $k$
and the codimension of the
immersion $V\to W$
equals that of $Z\to X$.
\end{lm}

\proof{
1. $h$ is $T^*X$-transversal
if and only if
$W\times_XT^*X\to T^*W$
is an injection.
This means that
$h^*\Omega^1_{X/k}\to \Omega^1_{W/k}$
is a locally split injection.

2. and 3.
Clear from the definition and 1.

4. 
It suffices to apply
Lemma \ref{lmbase}.2
to the closed conical subset
$h^*C\cap
{\rm Ker}(W\times_XT^*X
\to T^*W))$
of $W\times_XT^*X$.

5.
Condition (1) means that
$T^*_ZX\times_XW\to
T^*W$
is an injection.
In other words,
${\mathcal N}_{Z/X}\otimes_{{\mathcal O}_Z}
{\mathcal O}_V$
is locally a direct summand of
$\Omega^1_{W/k}
\otimes_{{\mathcal O}_W}{\mathcal O}_V$.
This is equivalent to condition (2).
\qed

}
\medskip

In the next lemma,
we regard
the closed subset
$h^*C$ as a closed subscheme.

\begin{lm}\label{lmhC}
{\rm (\cite[Lemma 1.2 (ii)]{SS})}
Let $h\colon W\to X$ be a 
morphism of smooth schemes over
$k$ and let $C\subset T^*X$ be
a closed conical subset.
Assume that
$h\colon W\to X$ is $C$-transversal.
Then, $W\times_XT^*X
\to T^*W$ is finite on $h^*C$.
\end{lm}

Under the assumption in Lemma \ref{lmhC},
we define
a closed conical subset
$h^\circ C\subset T^*W$
to be the image of $h^*C$.
Recall that $h^*C$ is the inverse image
of $C\subset T^*X$
by $T^*X\times_XW\to T^*X$.

\begin{lm}\label{lmhg}
Let $g\colon U\to X$ be a 
morphism of smooth schemes over
$k$ and let $C\subset T^*X$ be
a closed conical subset.
Assume that
$g$ is $C$-transversal.
Then 
for a morphism
$h\colon W\to U$ of
smooth schemes over $k$,
the following conditions
are equivalent:

{\rm (1)}
$g\circ h$ is $C$-transversal.

{\rm (2)}
$h$ is $g^\circ C$-transversal.
\end{lm}

\proof{
We consider the composition
$(g\circ h)^*C
\subset T^*X\times_XW
\to
T^*U\times_UW
\to T^*W$.
The assumption implies that
the intersection 
$(g\circ h)^*C
\cap{\rm Ker}( T^*X\times_XW
\to
T^*U\times_UW)$
is a subset of the 0-section.
Since 
$h^*(g^\circ C)
\subset T^*U\times_UW$
is the image of
$(g\circ h)^*C$,
the intersection 
$(g\circ h)^*C
\cap{\rm Ker}( T^*X\times_XW
\to
T^*W)$
is a subset of the 0-section
if and only
$h^*(g^\circ C)
\cap 
{\rm Ker}( T^*U\times_UW
\to
T^*W)$
is a subset of the 0-section.
\qed

}

\medskip
\begin{df}\label{dfgC}
Let $g\colon X'\to X$ be a 
morphism of smooth schemes over
$k$.

{\rm 1.}
Let $C'\subset T^*X'$
be a closed conical subset
and let $B'\subset X'$ be the base of
$C'$.
We say that $g$ is proper
on $B'$ if the restriction of $g$
on $B'$ regarded as a closed subscheme
of $X'$ is proper.

Assume that $g$ is proper
on $B'$.
We define
a closed conical subset
$g_\circ C'\subset T^*X$
to be the image by
$X'\times_XT^*X\to T^*X$
of the 
inverse image of $C'$
by the canonical morphism
$X'\times_XT^*X \to T^*X'$.

{\rm 2.}
Let
$h\colon W\to X$
be a morphism
of smooth schemes over $k$
and $V'\subset V=X'\times W$
be an open subscheme.
We say that $h$ and $g$
are transversal on $V'$
if $V'$ is smooth over $k$
and if the morphism
${\mathcal O}_W\otimes_{{\mathcal O}_X}^L
{\mathcal O}_Z
\to {\mathcal O}_V$
is an isomorphism on $V'$.
\end{df}

If $i\colon Z\to X$ is a closed
immersion of smooth schemes over
$k$ and if $C\subset T^*Z$
is a closed conical subset,
then  
$i_\circ C\subset T^*X$
is 
the inverse image of $C$
by the canonical surjection
$T^*X|_Z\to T^*Z$
regarded as a subset of $T^*X$.
Condition (2) in Lemma \ref{lmCtrans}.5
means that 
the morphism $h\colon W\to X$ 
is transversal to the immersion $i\colon Z\to X$.

\begin{lm}\label{lmiCtrans}
Let $g\colon X'\to X$ be a 
morphism of smooth schemes over
$k$ and let $C'\subset T^*X'$
be a closed conical subset.
Assume that $g$ is proper
on the base $B'\subset X'$ of
$C'$ and define $C=g_\circ C'$.
Let $h\colon W\to X$
be a morphism of smooth schemes
over $k$.
Assume that $h$ is $g_\circ C'$-transversal.

{\rm 1.}
There exists an open neighborhood 
$V'
\subset V=
X'
\times_XW$
of the inverse image
${\rm pr}_1^{-1}(B')$
such that
$g$ and $h$ are transversal
on $V'$.

{\rm 2.}
Let $V'\subset V$
be an open subset as in {\rm 1.}
Then,
the restriction of the projection
$h'\colon V'\to X'$
is $C'$-transversal
and $h^\circ g_\circ C'$
equals $g'_\circ h'^\circ C'$
for $g'\colon V'\to W$.
\end{lm}

\proof{
1.
By decomposing $h$
into $W\to X\times W\to X$,
it suffices to consider the case where $h$
is an immersion by Lemma \ref{lmhg}.
Since the immersion $h\colon W\to X$
is $g_\circ C'$-transversal,
the intersection
${\rm Ker}(T^*X\times_XX'\to T^*X')\times_{X'}V
\cap
{\rm Ker}(T^*X\times_XW\to T^*W)
\times_WV$
is a subset of the 0-section
on a neighborhood $V'$ of $
{\rm pr}_1^{-1}(B')\subset V$.
Namely, the restriction of
$(T^*X\times_XX'\to T^*X')\times_{X'}V'$
on 
${\rm Ker}(T^*X\times_XW\to T^*W)
\times_WV'$
is a locally split injection and
the assertion follows.

2.
Similarly as in the proof of 1,
we may assume that $h\colon W\to X$ is an immersion.
Then, further,
we have a commutative diagram
$$\begin{CD}
0@>>> T^*_WX\times_WV'
@>>> T^*X\times_XV'
@>>> T^*W\times_WV'
@>>>0
\\
@.@VVV@VVV@VVV@.\\
0@>>> T^*_{V'}X'
@>>> T^*X'|_{V'}
@>>> T^*V'
@>>>0
\end{CD}$$
of exact sequences
where the left vertical arrow is an isomorphism.
Hence the $C$-transversality of $h$
implies the $C'$-transversality of $h'\colon
V'\to X'$
and the equality 
$h^\circ g_\circ C'
=g'_\circ h'^\circ C'$ follows.
\qed

}

\begin{lm}\label{lmsC}
Let $X$ be a smooth scheme
over a field $k$
and let
$W\subset Z\subset X$
be closed subschemes
smooth over $k$.
Let $s\colon T^*Z\to T^*X|_Z$
be a linear section of the canonical
surjection
$T^*X|_Z\to T^*Z$.

{\rm 1.}
Locally on a neighborhood of $W\subset X$,
there exists a closed subscheme
$V\subset X$
smooth over $k$
such that $W\to Z\times_XV$
is an isomorphism,
${\rm codim}_XV=
{\rm codim}_ZW$
and $T^*_VX|_W\subset T^*X|_W$
equals the image
$s(T^*_WZ)$ by the section $s$.

{\rm 2.}
Let $C\subset T^*X|_Z$
be a closed conical subset.
Let $h\colon W\to Z$ be the closed immersion 
and 
let $h'\colon V\to X$
be a closed immersion as in {\rm 1.}
If $h$ is $s^{-1}(C)$-transversal,
then $h'$ is $C$-transversal on
a neighborhood of the image of $W$.
\end{lm}

\proof{
1.
Let ${\mathcal I}\subset {\mathcal O}_Z$
and ${\mathcal J}\subset {\mathcal O}_X$
be the ideal sheaves
defining $W\subset Z$
and $W\subset X$.
Then the section
$s\colon T^*Z\to T^*X|_Z$
induce a section
$T^*_WZ\to T^*_WX$
on sub vector bundles,
namely
${\mathcal I}/{\mathcal I}^2
\to {\mathcal J}/{\mathcal J}^2$.
Locally on $W$,
we take a basis of
$s({\mathcal I}/{\mathcal I}^2)
\subset {\mathcal J}/{\mathcal J}^2$.
By taking
its lifting to ${\mathcal J}$,
we
define a closed subscheme 
$V\subset X$.
Then, for the conormal sheaf,
we have an isomorphism
${\mathcal N}_{V/X}\otimes_{{\mathcal O}_V}
{\mathcal O}_W\to {\mathcal N}_{W/Z}$.
Hence after shrinking $V$ if necessary,
$V$ is smooth over $k$
and ${\rm codim}_XV=
{\rm codim}_ZW$.
The  isomorphism
${\mathcal N}_{V/X}\otimes_{{\mathcal O}_V}
{\mathcal O}_W\to {\mathcal N}_{W/Z}$
means the
equality $T^*_VX|_W=s(T^*_WZ)$.

2.
Since $V\cap Z=W$,
the equality $T^*_VX|_W=s(T^*_WZ)$
implies that the section
$s$ induces a bijection
$s^{-1}(C)\cap {\rm Ker}(T^*Z|_W\to T^*W)
\to
C\cap {\rm Ker}(T^*X|_V\to T^*V)$.
\qed

}


\section{$C$-acyclicity}

\begin{df}\label{dfCacyc}
{\rm (\cite[1.2]{SS})}
Let $X$
be a smooth scheme over a field $k$
and $C$ be a closed conical subset of 
the cotangent bundle $T^*X$.

{\rm 1.}
We say that
a morphism $f\colon X\to Y$
of smooth schemes over $k$
is $C$-acyclic
if the inverse image of
$C$ by the canonical morphism
$X\times_YT^*Y\to T^*X$
is a subset of the $0$-section.

{\rm 2.}
Let $h\colon W\to X$
and $f\colon W\to Y$
be morphisms
of smooth schemes over $k$.
We say that $(h,f)$ is $C$-acyclic
if 
$h\colon W\to X$
is $C$-transversal
and if $f\colon W\to Y$ is
$h^\circ C$-acyclic.

We say that $(h,f)$ is  universally $C$-acyclic
over $k$
if for every morphism
$g\colon Y'\to Y$ of smooth schemes over $k$
and 
for 
the commutative diagram
\begin{equation}
\xymatrix{
&W'\ar[r]^{f'}\ar[ld]_{h'}\ar[d]^{g'}&Y'\ar[d]^g\\
X&W\ar[l]_h\ar[r]^f &Y}
\label{eqCuniv}
\end{equation}
with cartesian square,
there exists a neighborhood $U'\subset W'$
of the inverse image
$h'^{-1}(B)$ of
the base $B$ of $C$
such that 
$g$ is transversal to $f$
on $U'$
and the pair $(h',f')$ is $C$-acyclic
on $U'$.
\end{df}

To avoid using the same terminology
for different notions,
we modified 
Beilinson's original terminology 
in \cite{SS}.

\begin{lm}\label{lmCacycsm}
Let $X$
be a smooth scheme over a field $k$
and $C$ be a closed conical subset of $T^*X$.

{\rm 1.}
Let $h\colon W\to X$
and $f\colon W\to Y$
be morphisms
of smooth schemes over $k$.
We identify $T^*(X\times Y)$
with $T^*X\times T^*Y$
and regard $C\times T^*Y$
as a closed conical subset of
$T^*(X\times Y)$.
Then, the following conditions are
equivalent:

{\rm (1)}
The pair $(h,f)$ is $C$-acyclic.

{\rm (2)}
The morphism $(h,f)\colon W\to X\times Y$
is $C\times T^*Y$-transversal.

\noindent
If these equivalent conditions
are satisfied,
then
the morphism $f\colon
W\to Y$ is smooth
on 
a neighborhood of
$h^{-1}(B)$.

{\rm 2.}
Let $f\colon X\to Y$
be a morphism of smooth schemes over $k$
and assume that
$C$ is the $0$-section $T^*_XX$.
Then 
$f\colon X\to Y$ is
$C$-acyclic
if and only if $f$ is smooth.
\end{lm}

\proof{
1.
(1)$\Rightarrow$(2):
Let $(a,b)$
be a point of the intersection
$\bigl(h^*C\times_W (T^*Y\times_YW)\bigr)
\cap 
{\rm Ker}((T^*X\times_XW)\times_W
(T^*Y\times_YW)
\to T^*W)$.
Then, $b\in T^*Y\times_YW$ 
is in the inverse image of
$-dh(a)\in h^\circ C$
and the $h^\circ C$-acyclicity
of $f$ implies that $b=0$.
Further the $C$-transversality of
$h$ implies $a=0$.

(2)$\Rightarrow$(1):
Since the injection
$T^*X\times_XW
\to
(T^*X\times_XW)\times_W
(T^*Y\times_YW)$
to the first factor
maps 
$h^*C
\cap 
{\rm Ker}(T^*X\times_XW
\to T^*W)$
to
$\bigl(h^*C\times_W (T^*Y\times_YW)\bigr)
\cap 
{\rm Ker}((T^*X\times_XW)\times_W
(T^*Y\times_YW)
\to T^*W)$,
the $C\times T^*Y$-transversality
of $(h,f)$ implies
the $C$-transversality
of $h$.
If $b\in T^*Y\times_YW$ 
is in the inverse image of
$dh(a)\in h^\circ C$ for $a\in h^* C$,
then 
$(-a,b)$ is contained in 
$\bigl(h^*C\times_W (T^*Y\times_YW)\bigr)
\cap 
{\rm Ker}((T^*X\times_YW)\times_W
(T^*Y\times_YW)
\to T^*W)$
and is $0$.

If $f\colon W\to Y$ is $h^\circ C$-acyclic,
then ${\rm Ker}(W\times_YT^*Y\to T^*W)|_{h^{-1}(B)}$
is a subset of the $0$-section.

2.
The morphism
$f$ is $T^*_XX$-acyclic
if and only if $f$
is $T^*Y$-transversal.
Hence the assertion
follows from
Lemma \ref{lmCtrans}.1.
\qed

}

\begin{lm}\label{lmCacyc}
Let $h\colon W\to X$
and $f\colon W\to Y$
be morphisms of smooth schemes over $k$
and $C$ be a closed conical subset of $T^*X$.

{\rm 1. (\cite[Example 1.2 (i)]{SS})}
The following conditions are
equivalent:

{\rm (1)}
$(h,f)$ is $T^*_XX$-acyclic.

{\rm (2)}
$f\colon W\to Y$ is smooth.

{\rm 2. (\cite[Example 1.2 (ii)]{SS})}
The following conditions are
equivalent:

{\rm (1)}
$(h,f)$ is $T^*X$-acyclic.

{\rm (2)}
$(h,f)\colon W\to X\times Y$ is smooth.
\end{lm}

\proof{
1.
Since $h$ is $T^*_XX$-transversal
by Lemma \ref{lmCtrans}.2
and $h^\circ T^*_XX=T^*_WW$,
condition (1) means that
$f\colon W\to Y$ is
$T^*_WW$-acyclic.
This is equivalent to (2)
by Lemma \ref{lmCacycsm}.2.

2.
By Lemma \ref{lmCacycsm}.1,
condition (1) means that
$(h,f)\colon W\to X\times Y$ is
$T^*X\times T^*Y$-transversal.
This is equivalent to (2)
by Lemma \ref{lmCtrans}.1.
\qed

}

\begin{lm}\label{lmhj}
Let $g\colon U\to X$ be a smooth
morphism of smooth schemes over
$k$ and let $C\subset T^*X$ be
a closed conical subset.
Let $g^\circ C$ be as defined before
Lemma {\rm \ref{lmhg}}.
Then 
for a pair of morphisms
$h\colon W\to U$ 
and $f\colon W\to Y$ of
smooth schemes over $k$,
the following conditions
are equivalent:

{\rm (1)}
$(g\circ h,f)$ is $C$-acyclic.

{\rm (2)}
$(h,f)$ is $g^\circ C$-acyclic.
\end{lm}

\proof{
By Lemma \ref{lmCacycsm}.1,
the conditions are equivalent to
the conditions:
{\rm (1)}
$(g\times 1_Y)\circ (h,f)$ is $C\times T^*Y$-transversal,
and 
{\rm (2)}
$(h,f)$ is $(g\times 1_Y)^\circ (C\times T^*Y)$-transversal,
respectively.
Hence they are equivalent by
Lemma \ref{lmhg}.
\qed

}

\begin{lm}\label{lmiCacyc}
Let $g\colon X'\to X$ be a 
morphism of smooth schemes over
$k$ and let $C'\subset T^*X'$
be a closed conical subset.
Assume that $g$ is proper
on the base $B'\subset X'$ of
$C'$
and define
$g_\circ C'$ as in Definition
{\rm \ref{dfgC}}.
Let $f\colon X\to Y$
be a morphism of smooth schemes
over $k$.
Then, the following conditions
are equivalent:

{\rm (1)}
$f$ is $g_\circ C'$-acyclic.

{\rm (2)}
$f\circ g$ is $C'$-acyclic.
\end{lm}

\proof{
Since $g_\circ C'\subset T^*X$ is defined
as the image of
the inverse image
of $C'$ by
$T^*X\times_XX'\to T^*X'$,
either of the conditions means that
the inverse image of $C'$ by
$T^*Y\times_YX'\to T^*X'$
is a subset of the $0$-section.
\qed

}

\begin{lm}\label{lmCacych}
Let $$\xymatrix{
X&W\ar[l]_h\ar[r]^f\ar[rd]_{f'}&Y\ar[d]^g\\
&&Y'}$$ 
be a commutative diagram
of smooth schemes of finite type over $k$.
Let $C\subset T^*X$ be a closed conical subset.
Assume that $g$ is smooth.
If $(h,f)$ is $C$-acyclic, then
$(h,f')$ is $C$-acyclic.
\end{lm}

\proof{
Since the inverse image of
$h^\circ C$ by
$T^*Y\times_YW\to T^*W$
is a subset of the $0$-section
and since
$T^*Y'\times_{Y'}Y\to T^*Y$
is an injection,
the inverse image of
$h^\circ C$ by
$T^*Y'\times_{Y'}W\to T^*W$
is a subset of the $0$-section.
\qed

}

\begin{lm}\label{lmCacycpr}
Let $g\colon X'\to X$ be a morphism
of smooth schemes over $k$.
Let $C'\subset T^*X'$ be a closed
conical subset and assume that $g$ is
proper on the base $B'$ of  $C'$.
Let $C=g_\circ C'$.
Let $$
\xymatrix{
X'\ar[d]_g&W'\ar[l]_{h'}\ar[rd]^{f'}\ar[d]_{g'}&\\
X&W\ar[l]_h\ar[r]^f &Y}
$$ be a commutative
diagram of separated morphisms
of schemes of finite
type over $k$ with cartesian square.
Assume that $W$ and $Y$ are
smooth over $k$ and
that
$(h,f)$ is $C$-acyclic.

Then, there exists a open neighborhood $U'
\subset W'$ of $h'^{-1}(B')$ smooth over $k$
on which
$(h',f')$ is $C'$-acyclic.
\end{lm}

\proof{
By Lemma \ref{lmiCtrans},
there exists an open neighborhood $U'
\subset W'$ of $h'^{-1}(B')$ 
such that $g$ and $h$ are transversal
on $U'$
and $h^\circ g_\circ C'=g'|_{U'\circ}
h'^\circ C'$.
Since $f$ is $h^\circ g_\circ C'$-acyclic,
$f'$ is $h'^\circ C'$-acyclic by
Lemma \ref{lmiCacyc}.
\qed

}

\begin{lm}\label{lmuniCacy}
Let $h\colon W\to X$ and
$f\colon W\to Y$ be morphisms of
smooth schemes of finite type
over $k$.
Let $C\subset T^*X$ be a closed conical
subset.
If $(h,f)$ is $C$-acyclic,
then $(h,f)$ is universally $C$-acyclic.
\end{lm}

\proof{
Let $g\colon Y'\to Y$
be a morphism of smooth schemes
over $k$ and
consider the commutative diagram (\ref{eqCuniv}).
By Lemma \ref{lmCacycsm}.1,
$f\colon W\to Y$ is smooth
on a neighborhood $U\subset W$ of $h^{-1}(B)$.
Hence 
the morphism $g$ is transversal to $f$ on 
the inverse image $U'\subset W'$ of $U$.
To show that $(h',f')$ is $C$-acyclic on $U'$,
it suffices to show
that 
\begin{align*}
&{\rm Ker}((T^*X\times_XW)
\times_W
(T^*Y\times_YW)
\to T^*W)\times_WU'\\
&\to 
{\rm Ker}((T^*X\times_XU')
\times_{U'}
(T^*Y'\times_{Y'}U')
\to T^*U')
\end{align*}
is an isomorphism
further by Lemma \ref{lmCacycsm}.1.
By the factorization
$Y'\to Y\times Y'\to Y$,
it suffices to prove the cases
where $Y'\to Y$ is smooth
and $Y'\to Y$ is an immersion
respectively.
Since the morphism $f\colon W\to Y$
is transversal to $g\colon Y'\to Y$ 
in the second case, the assertion follows.
\qed

}

\section{Micro support}

Let $f\colon X\to S$ be a morphism
of schemes. Let
$x\to X$ and $t\to S$ be geometric points
and let $S_{(s)}$ be the strict localization
of $S$ at the image $s=f(x)$ of $x$.
A specialization $s\gets t$ means
a lifting of $t\to S$ to $t\to S_{(s)}$.

In this article, we assume that
a prime number $\ell$ is invertible
on the considered scheme $X$ 
and let $\Lambda$ denote a finite field
of characteristic $\ell$.
By abuse of terminology,
an object ${\mathcal F}$ of
$D^+(X,\Lambda)$ will be called
a sheaf and it is called
constructible if 
the cohomology sheaf
${\mathcal H}^q{\mathcal F}$
is constructible for every $q
\in {\mathbf Z}$
and 
is 0 except for 
finitely many $q$.
We say that a sheaf ${\mathcal F}$ is locally constant if
every cohomology sheaf
${\mathcal H}^q{\mathcal F}$ is locally constant.

\begin{df}\label{dfla}
Let $f\colon X\to S$ be a morphism
of schemes
and ${\mathcal F}$ be a sheaf on $X$.
We say that $f$ is locally acyclic relatively to ${\mathcal F}$
or ${\mathcal F}$-acyclic for short
if for each geometric points
$x\to X$ and $t\to S$
and each specialization $s=f(x)\gets t$,
the canonical morphism
${\mathcal F}_x\to 
R\Gamma(X_{(x)}\times_{S_{(s)}}t,
{\mathcal F})$
is an isomorphism.

We say that $f$ is universally 
${\mathcal F}$-acyclic,
if for every morphism $S'\to S$,
the base change of $f$ is locally acyclic
relatively to the pull-back of ${\mathcal F}$.
\end{df}

The local acyclicity
is an \'etale local property on $X$
by definition. 

For geometric points $s,t$ of $S$
and a specialization $t\to S_{(s)}$,
let $i\colon X_s\to X\times_SS_{(s)}$
and $j\colon X_t\to X\times_SS_{(s)}$
denote the canonical morphisms.
Then, the local acyclity is equivalent to
the condition
that the canonical morphism
$i^*{\mathcal F}_{X\times_SS_{(s)}}\to i^*j_*j^*{\mathcal F}_{X\times_SS_{(s)}}$
is an isomorphism for each $s,t$ and $s\gets t$.
In fact,
the morphism
${\mathcal F}_x\to 
R\Gamma(X_{(x)}\times_{S_{(s)}}t,
{\mathcal F})$
in Definition \ref{dfla} is
the stalk at $x$ of
$i^*{\mathcal F}_{X\times_SS_{(s)}}\to i^*j_*j^*{\mathcal F}_{X\times_SS_{(s)}}$.

%

\begin{thm}\label{thmlcc}
Let $X$ be a scheme
and ${\mathcal F}$ be a sheaf on $X$.

{\rm 1}.
We consider the following conditions:

{\rm (1)}
${\mathcal F}$ is locally constant.

{\rm (2)}
Every smooth morphism
$f\colon X\to Y$ 
is ${\mathcal F}$-acyclic.

{\rm (3)}
The identity $1_X\colon X\to X$
is ${\mathcal F}$-acyclic.

We have {\rm (1)}$\Rightarrow${\rm (2)}$\Rightarrow${\rm (3)}.
If $X$ is a locally noetherian scheme
and if, for every geometric point $x$
of $X$, the stalk ${\mathcal F}_x$ is constructible,
then the three conditions are equivalent.

{\rm 2.
(\cite[Corollaire 2.16]{TF})}
If $X$ is a scheme over $k$,
then the morphism
$X\to {\rm Spec}\, k$
is universally ${\mathcal F}$-acyclic.
\end{thm}

\proof{
1.
The implication
(1)$\Rightarrow$(2)
is 
the local acyclicity of smooth morphisms
\cite[Th\'eor\`eme 2.1]{AL}
and
(2)$\Rightarrow$(3)
is clear.
The implication
(3)$\Rightarrow$(1)
is \cite[Proposition 2.13 (i)]{Artin}.
\qed

}

%
%
%
%
%
%

\begin{lm}\label{lmFacycg}
%
Let
$$\xymatrix{
X'\ar[d]_g
\ar[rd]^{f'}&\\
X
\ar[r]^f &Y}
$$
be a commutative diagram
of schemes of finite type over $k$.
Let ${\mathcal F}'$ be
a constructible sheaf on $X'$
and $B'\subset X'$ be 
its 
support.
Assume that $g$ is proper 
on $B'$
and let ${\mathcal F}=g_*{\mathcal F}'$.
We consider the following conditions:

{\rm (1)} $f$ is ${\mathcal F}$-acyclic over $S$.

{\rm (2)} 
$f'$ is ${\mathcal F}'$-acyclic over $S$.

We have the implication {\rm (2)}$\Rightarrow${\rm (1)}.
If $g$ is finite
on $B$, 
then these conditions are equivalent.
\end{lm}

\proof{
Let $s,t$ be geometric points  of $Y$
and $t\to Y_{(s)}$ be a specialization
and consider the cartesian diagram
$$\begin{CD}
X'_s@>{i'}>> X'\times_YY_{(s)}
@<{j'}<<X'_t\\
@V{g'_s}VV@VV{g'_{Y_{(s)}}}V@VV{g'_t}V\\
X_s@>i>> X\times_YY_{(s)}
@<j<<X_t
\end{CD}$$
where the vertical arrows are induced
by $g'$.
Since $g$ is assumed proper on the support $B'$
of ${\mathcal F}'$,
the vertical arrows in the commutative
diagram
$$\begin{CD}
g'_{s*}i'^*{\mathcal F}'_{X'\times_YY_{(s)}}
@>>>
g'_{s*}i'^*j'_*j'^*{\mathcal F}'_{X'\times_YY_{(s)}}\\
@AAA@AAA\\
i^*{\mathcal F}_{X\times_YY_{(s)}}
@>>>
i^*j_*j^*{\mathcal F}_{X\times_YY_{(s)}}
\end{CD}$$
are isomorphisms.
If {\rm (2)} is satisfied,
then the upper horizontal arrow is an
isomorphism.
Hence the lower horizontal arrow is an
isomorphism and {\rm (1)} holds.
Conversely, if $g$ is finite on $B$
and if the lower horizontal arrow is an
isomorphism,
the upper horizontal arrow is an isomorphism
before taking $g'_{s*}$
and {\rm (1)} holds.
\qed

}

\begin{df}
Let $X$ be a smooth scheme
over a field $k$.
Let ${\mathcal F}$ be a sheaf on $X$
and let $C\subset T^*X$ be a closed conical subset.
We say that ${\mathcal F}$ is micro supported on $C$,
if for every $C$-acyclic pair $(h,f)$ of morphisms
$h\colon W\to X$ and $f\colon W\to Y$ of smooth schemes
over $k$,
the morphism $f$ is $h^*{\mathcal F}$-acyclic.
\end{df}

If ${\mathcal F}$ is micro supported on $C\subset C'$,
then
${\mathcal F}$ is micro supported on $C'$
by Lemma \ref{lmCacycsm}.1
and Lemma \ref{lmCtrans}.3
(cf.~\cite[1.3]{SS}).

\begin{lm}\label{lmms}
Let $X$ be a smooth scheme
over a field $k$
and let ${\mathcal F}$ be a sheaf on $X$.

{\rm 1. (cf.~\cite[Lemma 2.1 (iii)]{SS})}
Assume that for every geometric point $x$
of $X$, the stalk ${\mathcal F}_x$ is constructible.
Then the following conditions are equivalent:

{\rm (1)} ${\mathcal F}$ is locally constant.

{\rm (2)} ${\mathcal F}$ is micro supported on
the $0$-section $T^*_XX$.

{\rm 2. (cf.~\cite[Lemma 1.3]{SS})}
Every sheaf ${\mathcal F}$ is
micro supported on $T^*X$.

{\rm 3.}
Let $C\subset
T^*X$ be a closed conical subset
on which ${\mathcal F}$ is micro supported
and $U\subset X$ be
an open subset.
If $C_U$ is empty,
then the restriction ${\mathcal F}_U$ is $0$.

Conversely, ${\mathcal F}=0$
is micro supported on $\varnothing$.

{\rm 4.}
Let $h\colon W\to X$ and
$f\colon W\to Y$ be morphisms 
of smooth schemes over $k$.
Let $C\subset
T^*X$ be a closed conical subset
on which ${\mathcal F}$ is micro supported.
If $(h,f)$ is $C$-acyclic,
then
$f$ is universally $h^*{\mathcal F}$-acyclic.
\end{lm}

Note that in the convention of \cite{SS},
a sheaf is assumed constructible.

\proof{
1.
By Lemma \ref{lmCacyc}.1,
a pair $(h,f)$ of morphisms
$h\colon W\to X$ and $f\colon W\to Y$ of smooth schemes
over $k$,
is $T^*_XX$-acyclic
if and only if $f$ is smooth.

(1)$\Rightarrow$(2):
If $f$ is smooth and ${\mathcal F}$ is locally constant,
then $f$ is $h^*{\mathcal F}$-acyclic by
Theorem \ref{thmlcc}.1 (1)$\Rightarrow$(2).

(2)$\Rightarrow$(1):
Since the pair
$(1_X,1_X)$ is $T^*_XX$-acyclic,
the condition (2) implies that
the identity $1_X$ is ${\mathcal F}$-acyclic.
By Theorem \ref{thmlcc}.1 (3)$\Rightarrow$(1),
this implies that
${\mathcal F}$ is locally constant.

2.
Suppose that a pair $(h,f)$
of morphisms $h\colon W\to X$
and $f\colon W\to Y$ of
smooth schemes over $k$,
is $T^*X$-acyclic.
Then $(h,f)\colon W\to X\times Y$
is smooth by Lemma \ref{lmCacyc}.2.
Since the question is \'etale local 
on $W$,
we may assume that $W={\mathbf A}^n_{X\times Y}$.
Then, we have a cartesian diagram
$$\begin{CD}
X@<<< X\times Y@<<< W\\
@VVV@VVV@VVV\\
{\rm Spec}\, k@<<< Y@<<< {\mathbf A}^n_Y.
\end{CD}$$
By Theorem \ref{thmlcc}.2,
the morphism
$W\to {\mathbf A}^n_Y$
is ${\mathcal F}_W$-acyclic
for the pull-back ${\mathcal F}_W$ of ${\mathcal F}$.
Since ${\mathbf A}^n_Y\to Y$
is smooth,
the composition
$W\to Y$ is also
${\mathcal F}_W$-acyclic
by \cite[Corollaire 2.7]{App}.

3.
Since the pair $(1,0)$
of the identity $1\colon U\to U$
and the $0$-mapping
$U\to {\mathbf A}^1$
is $C$-acyclic,
the constant morphism
$0\colon U\to {\mathbf A}^1$
is ${\mathcal F}$-acyclic.
Hence we have ${\mathcal F}=0$.

If ${\mathcal F}=0$,
every pair $(h,f)$ of morphisms
$h\colon W\to X$ and
$f\colon W\to Y$ of smooth schemes
over $k$ is ${\mathcal F}$-acyclic.

4.
Since $(h,f)$ is universally $C$-acyclic
by Lemma \ref{lmuniCacy},
the assertion follows.
\qed

}

\begin{lm}\label{lmmsUi}
{\rm (cf.~\cite[Lemma 2.1 (iv)]{SS})}
Let $X$ be a smooth scheme
over a field $k$
and ${\mathcal F}'\to {\mathcal F}\to {\mathcal F}''\to$
be a distinguished triangle of sheaves
on $X$.
Suppose that
${\mathcal F}$ and ${\mathcal F}'$ are micro supported on
$C$ and on $C'$ respectively.
Then
${\mathcal F}''$ is micro supported on $C''
=C\cup C'$.
\end{lm}

\proof{
Let $h\colon W\to X$ and
$f\colon W\to Y$ be morphisms
of smooth schemes over $k$
such that the pair $(h,f)$ is $C''$-acyclic.
Then, since $C,C'\subset C''$,
the pair $(h,f)$ is $C$-acyclic
and $C'$-acyclic.
Hence $f$ is $h^*{\mathcal F}$-acyclic
and $h^*{\mathcal F}'$-acyclic.
By the distinguished triangle
$ {\mathcal F}'\to {\mathcal F}\to {\mathcal F}''\to$,
$f$ is $h^*{\mathcal F}''$-acyclic.
\qed

}

\begin{lm}\label{lmmsloc}
Let $X$ be a smooth scheme
over a field $k$
and let ${\mathcal F}$ be a sheaf on $X$.

{\rm 1. (cf.~\cite[Lemma 2.2 (i)]{SS})}
Assume that ${\mathcal F}$ is micro supported on $C$.
If a morphism $g\colon V\to X$ of smooth schemes over $k$ is $C$-transversal,
then the pull-back
$g^*{\mathcal F}$ is micro supported on $g^\circ C$.

{\rm 2.}
Let $(j_i\colon U_i\to X)_{i\in I}$
be an \'etale covering
and let $C_i\subset T^*U_i$
be a family of closed conical subsets
such that the pull-backs
${\mathcal F}_{U_i}$ are micro supported
on $C_i$.
We identify $T^*U_i$
with $T^*X\times_XU_i$
and consider the image of
$C_i\subset T^*U_i=T^*X\times_XU_i
$ as a subset $
{\rm Im}(C_i) \subset T^*X$.
Then, 
${\mathcal F}$ is micro supported
on the closure $C=\overline{
\bigcup_{i\in I}
{\rm Im}(C_i)}$
of the union of the images.

{\rm 3.}
Let $U\subset X$ be an open subset
and $Z=X\sm U$ be the complement.
Assume that ${\mathcal F}$ is
micro supported on $C\subset T^*X$
and that 
the restriction ${\mathcal F}_U$ is
micro supported on $C'\subset T^*U$.
Then, ${\mathcal F}$ is
micro supported on the union
$\overline {C'}\cup C|_Z\subset T^*X$.
\end{lm}

\proof{
1.
Let $h\colon W\to V$
and $f\colon W\to Y$
be morphisms of smooth schemes over $k$
such that the pair $(h,f)$
is $g^\circ C$-acyclic.
Then, $(gh,f)$ is $C$-acyclic
by Lemma \ref{lmhj} (2)$\Rightarrow$(1).
Hence
 $f$ is locally acyclic relatively to $(gh)^*{\mathcal F}=h^*(g^*{\mathcal F})$.

2.
Let $h\colon W\to X$ and
$f\colon W\to Y$ be morphisms
of smooth schemes over $k$
such that the pair $(h,f)$ is $C$-acyclic.
Then, for every $i\in I$,
the pair of 
$h_i\colon W_i=W\times_XU_i\to U_i$
and $f_i\colon W_i\to Y$
is $j_i^\circ C$-acyclic
by Lemma \ref{lmhj} (2)$\Rightarrow$(1).
Further it is $C_i$-acyclic
since $C_i\subset j_i^\circ C$
by Lemma \ref{lmCtrans}.3.
Hence $f_i$ is locally acyclic
relatively to
$h_i^*({\mathcal F}_{U_i})=(h^*{\mathcal F})_{W_i}$
for every $i\in I$.
Since $(W_i\to W)_{i\in I}$
is an \'etale covering,
$f$ is $h^*{\mathcal F}$-acyclic.

{\rm 3.}
Let $h\colon W\to X$ and $f\colon W\to Y$
be morphisms of smooth schemes over $k$
and assume that $(h,f)$ is
$\overline {C'}\cup C|_Z$-acyclic.
Then, $(h,f)$ is $\overline {C'}$-acyclic
on $U$ and
is $C$-acyclic on
a neighborhood $V$ of $Z$.
Hence $f$ is $h^*{\mathcal F}$-acyclic
on $X=U\cup V$.
\qed

}

\begin{lm}\label{lmCA}
{\rm (cf.~\cite[Lemma 2.5]{SS})}
Let $X$ be a smooth scheme
over a field $k$ and let
${\mathcal F}$ be a sheaf on $X$.
Let $i\colon Z\to X$ be a closed immersion.
Assume that
${\mathcal F}
=i_*{\mathcal F}_Z$ is
supported on $Z$
and 
is micro supported on 
a closed conical subset $C\subset T^*X$.

{\rm 1. (cf.~\cite[Lemma 2.3]{SS})}
The sheaf
${\mathcal F}$ is micro supported on 
$C|_Z$.

{\rm 2.}
Assume that $Z$ is smooth
over $k$
and that
$C=C|_Z$ is a subset of
$T^*X|_Z$.
Let $s\colon T^*Z\to T^*X|_Z$
be a section of the surjection $T^*X|_Z\to T^*Z$.
Then, 
${\mathcal F}_Z$ is micro supported on
$s^{-1}(C)$.
Define $C_Z\subset T^*Z$ 
to be the closure of
the image of $C$
by the surjection
$T^*X|_Z\to T^*Z$.
Then, 
${\mathcal F}_Z$ is micro supported on
$C_Z$.

{\rm 3. (cf.~\cite[Lemma 2.2 (ii)]{SS})}
Assume that ${\mathcal F}_Z$
is micro supported on
a closed conical subset $C_Z\subset T^*Z$.
Then ${\mathcal F}$ is micro supported on $i_\circ C_Z$.
\end{lm}

\proof{
1.
Let $h\colon W\to X$ and
$f\colon W\to Y$ be morphisms
of smooth schemes over $k$
such that the pair $(h,f)$ is $C|_Z$-acyclic.
Then, by Lemma \ref{lmCtrans}.4,
there exists an open neighborhood 
$U$ of $Z$ such that the pair
of $h_U\colon W\times_XU\to X$
and $f_U\colon W\times_XU\to Y$
is $C$-acyclic.
Hence $f_U$ is
$(h^*{\mathcal F})_{W\times_XU}$-acyclic.
For the complement $V=X\sm Z$, 
the restriction $(h^*{\mathcal F})_{W\times_XV}$
is $0$ 
since ${\mathcal F}
=i_*{\mathcal F}_Z$.
Hence
$f_V\colon W\times_XV\to Y$
is $(h^*{\mathcal F})_{W\times_XV}$-acyclic.
Since
$U\cup V=X$,
$f$ is ${\mathcal F}$-acyclic.

2.
First, we show that
${\mathcal F}_Z$ is micro supported on 
$s^{-1}(C)$.
Let $h\colon W\to Z$ and
$f\colon W\to Y$ be morphisms
of smooth schemes over $k$
such that the pair $(h,f)$ is 
$s^{-1}(C)$-acyclic.
We show that we may assume $h$ is
an immersion.
We consider the commutative
diagram
$$\xymatrix{
X\times W\ar[d]_{{\rm pr}_1}
&
Z\times W\ar[l]\ar[d]_{{\rm pr}_1}
&&\\
X&Z\ar[l]
&W\ar[r]\ar[lu]\ar[l]
&Y.}
$$
For the section
$\widetilde s=s\times 1\colon
T^*Z\times T^*W=T^*(Z\times W)
\to T^*X|_Z\times T^*W=T^*(X\times W)|_{Z\times W}$
extending $s$,
we have
${\rm pr}_1^\circ(s^{-1}(C))
=\widetilde s^{-1}({\rm pr}_1^\circ C)$.
Since the projection
$Z\times W\to Z$ is smooth,
we may replace $h\colon W\to Z$,
$X\supset Z$ and $C$ 
by $(h,1)\colon W\to Z\times W$,
$X\times W
\supset Z\times W$
and ${\rm pr}_1^\circ C$,
by Lemma \ref{lmmsloc}.1.
Thus
we may assume that $h$ is an immersion.

Since the assertion is local on $X$,
by Lemma \ref{lmsC}.1,
we may
extend the immersion 
$(h,f)\colon W\to Z\times Y$
to an immersion 
$(h',f')\colon V\to X\times Y$ 
transversal to the immersion 
$Z\times Y\to X\times Y$
such that 
the restriction
$T^*_V(X\times Y)|_W$ 
of the conormal bundle is 
the image of 
$T^*_W(Z\times Y)$ by
the section 
$s|_W\times 1\colon
T^*Z|_W\times_W (T^*Y\times_YW)
\to T^*X|_W\times_W (T^*Y\times_YW)$.

We consider the cartesian diagram
$$\begin{CD}
W@>>> Z\times Y\\
@VVV@VVV\\
V@>>> X\times Y.\end{CD}$$
Since
$T^*_V(X\times Y)|_W$ is 
the image of 
$T^*_W(Z\times Y)$ by
the section,
the $s^{-1}(C)$-acyclicity of 
$(h,f)\colon W\to Z\times Y$
implies the $C$-acyclicity of 
$(h',f')\colon V\to X\times Y$
by Lemma \ref{lmCacycsm}.1
and Lemma \ref{lmsC}.2,
after shrinking $V$ if necessary.
Since ${\mathcal F}$ is micro supported on $C$,
this implies that 
$f'\colon V\to Y$ is 
$h'^*{\mathcal F}$-acyclic.
This means
that $f\colon W\to Y$ is 
$h^*{\mathcal F}_Z$-acyclic
since $h'^*{\mathcal F}$
is the direct image of
$h^*{\mathcal F}_Z$ by the closed immersion
$W\to V$.
Hence ${\mathcal F}_Z$ is micro supported on 
$s^{-1}(C)$.

We show that
${\mathcal F}_Z$ is micro supported on 
$C_Z$.
Since the assertion is local on $Z$,
we may assume that there exists
a section $s\colon T^*Z\to T^*X|_Z$.
Then, since
$s^{-1}(C)\subset C_Z$,
the assertion follows.

3.
Let $h\colon W\to X$ and
$f\colon W\to Y$ be morphisms
of smooth schemes over $k$
such that the pair $(h,f)$ is 
$i_\circ C$-acyclic.
By Lemma \ref{lmms}.3,
the support of ${\mathcal F}$
is a subset of the base $B$ of
$C$.
By Lemma \ref{lmiCtrans}
and Lemma \ref{lmiCacyc},
there exists an open neighborhood
$V'\subset V=Z\times_XW$ of the
inverse image of $B$
such that $V'$ is smooth over $k$,
that the pair of
$h'\colon V'\to Z$ 
and the composition 
$f'\colon V'\to W\to Y$ is
$C$-acyclic.
Hence $f'$ is $h'^*{\mathcal F}$-acyclic.
Since $V'$ contains
the inverse image of
the support of ${\mathcal F}$,
the morphism $f$ is $h^*i_*{\mathcal F}$-acyclic.
\qed

}

\begin{lm}\label{lmCbarms}
Let $g\colon X'\to X$ be a morphism
of smooth schemes over $k$.
Assume that $g$ is proper on the base of 
a closed conical subset $C'
\subset T^*X'$.
Let ${\mathcal F}'$ be a constructible sheaf on $X'$.
If ${\mathcal F}'$ is micro supported on $C'$,
then $Rg_*{\mathcal F}'$ is micro supported
on $C=g_\circ C'$.
\end{lm}

\proof{
Let $h\colon W\to X$ and $f\colon W\to Y$
be morphisms of smooth schemes over $k$
such that $(h,f)$ is $C$-acyclic.
Let $$
\begin{CD}
X'@<{h'}<<W'\\
@VgVV@VV{g'}V\\
X@<h<<W
\end{CD}$$
be the cartesian diagram.
Then, by Lemma \ref{lmCacycpr},
there exists an open regular neighborhood
$U'\subset W'=W\times_XX'$
of the inverse image of $B'$
where $h$ is transversal to $g$ and
$(h',fg')$ is $C'$-acyclic.
Hence the assertion follows from
Lemma \ref{lmFacycg}.
\qed

}

\section{Singular support and Radon transform}

\begin{df}
{\rm (\cite[1.3]{SS})}
Let $X$ be a smooth scheme
over a field $k$
and ${\mathcal F}$ be a sheaf on $X$.
We say that a closed conical
subset $C\subset T^*X$
is the singular support of ${\mathcal F}$
if for every closed conical subset
$C'\subset T^*X$,
the inclusion $C\subset C'$
is equivalent to the condition
that ${\mathcal F}$ is micro supported on $C'$.
If the singular support of ${\mathcal F}$
exists,
it is denoted by
$SS{\mathcal F}$.
\end{df}

Following Beilinson's argument in \cite{SS},
we will prove the existence of
singular support
for any sheaf without assuming
constructiblity
by reducing to the case
where $X$ is the projective space
${\mathbf P}^n$
and 
using the Radon transform
in the case $X={\mathbf P}^n$.
For the sake of convenience and
completeness, we give the detail.

\begin{lm}\label{lmSSloc}
Let $X$ be a 
smooth scheme over $k$ and
let ${\mathcal F}$ be a sheaf on $X$.

{\rm 1.}
Let $U\subset X$ be an open subscheme.
Assume that $C\subset T^*X$ is the singular support
of ${\mathcal F}$.
Then, $C|_U$ is the singular support
of ${\mathcal F}|_U$.

{\rm 2.}
Let $(U_i)$ be an open covering of
$X$ and $C_i
\subset T^*U_i$ be the singular support
of ${\mathcal F}|_{U_i}$.
Then,
$C=\bigcup_iC_i
\subset T^*X$ is the singular support
of ${\mathcal F}$.
\end{lm}

\proof{
1.
The restriction ${\mathcal F}|_U$
is micro supported on
$C|_U$ by Lemma \ref{lmmsloc}.1.
We show that $C|_U$ is the smallest.
Let $C'\subset T^*U$ be a closed conical
subset on which ${\mathcal F}|_U$
is micro supported
and $Z=X\sm U$ be the complement.
Then ${\mathcal F}$ is micro supported
on $\overline {C'}\cup C|_Z$
by Lemma \ref{lmmsloc}.3.
Since $C$ is the smallest,
we have
$C\subset 
\overline {C'}\cup C|_Z$
and
$C|_U\subset 
(\overline {C'}\cup C|_Z)|_U=C'$.

2.
For every $i,j$,
the restrictions $C_i|_{U_i\cap U_j}$ and
$C_j|_{U_i\cap U_j}$
are the singular support of 
${\mathcal F}|_{U_i\cap U_j}$
and are the same by 1.
Hence the union 
$C=\bigcup_iC_i$ is a closed conical subset of $T^*X$.
Since $C_i=C|_{U_i}$,
${\mathcal F}$ is micro supported on $C$
by Lemma \ref{lmmsloc}.1.

We show that $C$ is the smallest.
Let $C'\subset
T^*X$ be a closed conical subset 
on which
${\mathcal F}$ is micro supported.
Then, for each $i$,
we have $C_i\subset C'|_{U_i}$.
Hence we have $C\subset C'$.
\qed
}

\begin{pr}\label{prSSXP}
Let $i\colon X\to P$ be a closed immersion
of smooth schemes of
over $k$ and
let ${\mathcal F}$ be a constructible
sheaf on $X$.
Let $C_P\subset T^*P$ 
be a closed conical subset
and
assume that $C_P$ is the singular support
of $i_*{\mathcal F}$.
Then the following holds.

{\rm 1.}
$C_P$ is a subset of
$T^*P|_X$.

{\rm 2.}
Define $C\subset T^*X$ 
to be the closure of
the image of $C_P$
by the surjection
$T^*P|_X\to T^*X$.
Then $C$ is the singular support $SS{\mathcal F}$
and we have $C_P=i_\circ C$.
\end{pr}


\proof{
1.
Let $U=P\sm X$ be the complement.
Since ${\mathcal F}|_U=0$
is micro supported on $\varnothing$
by Lemma \ref{lmms}.3,
${\mathcal F}$ is micro supported on
$C_P|_X\subset T^*P$ 
by Lemma \ref{lmmsloc}.3.
Since $C_P$ is the smallest,
we have $C_P=C_P|_X
\subset T^*P|_X$.

2.
By Lemma \ref{lmCA}.2,
${\mathcal F}$ is micro supported on $C$.
We show that $C$ is the smallest.
Assume that ${\mathcal F}$ is micro supported on 
a closed conical subset
$C'\subset T^*X$.
Then, since $i_*{\mathcal F}$ is micro supported on $i_\circ C'$
by  Lemma \ref{lmCA}.3,
we have $C_P\subset i_\circ C'$.
By taking the closure of the image 
by the surjection
$T^*P|_X\to T^*X$,
we obtain $C\subset C'$.

Since $i_*{\mathcal F}$ is micro supported on $i_\circ C$,
we have $C_P\subset i_\circ C$.
We show
the other inclusion
$i_\circ C\subset C_P$.
Since the assertion is local,
we may assume that
there exists 
a section
$s\colon T^*X\to T^*P|_X$
of the surjection
$T^*P|_X\to T^*X$.
Since ${\mathcal F}$
is micro supported on
$s^{-1}(C_P)$ 
by Lemma \ref{lmCA}.2,
we have
$C\subset s^{-1}(C_P)$.
Hence we have
$i_\circ C\subset C_P$.
\qed

}


\medskip

We recall the Radon transform and the Legendre transform.
Let $V={\mathbf A}^{n+1}$
and let  ${\mathbf P}={\mathbf P}(V)$
be the projective space of dimension $n$
parametrizing lines in $V$.
The dual projective space ${\mathbf P}^\vee=
{\mathbf P}(V^\vee)$
is the moduli space of hyperplanes
in ${\mathbf P}$.

By the exact sequence
$0\to \Omega^1_{{\mathbf P}/k}(1)
\to {\mathcal O}_{\mathbf P}\otimes V^\vee
\to {\mathcal O}_{\mathbf P}(1)\to 0$
of locally free ${\mathcal O}_{\mathbf P}$-modules,
we define a closed subscheme
$Q={\mathbf P}(T^*{\mathbf P})
\subset 
{\mathbf P}\times {\mathbf  P}^\vee$
of codimension 1.
This equals the universal family of hyperplanes
since it is defined by the tautological section
$\Gamma({\mathbf P}\times {\mathbf  P}^\vee,
{\mathcal O}(1)\boxtimes{\mathcal O}(1))
=V^\vee\otimes V$ corresponding
to the identity $1\in {\rm End}(V)$.
Let $q\colon Q\to {\mathbf P}$
and $q^\vee\colon Q\to {\mathbf P}^\vee$
be the restrictions of the projections
${\mathbf P}\times {\mathbf  P}^\vee
\to {\mathbf P}$
and
${\mathbf P}\times {\mathbf  P}^\vee
\to {\mathbf P}^\vee$.
By symmetry, $Q\subset {\mathbf P}\times {\mathbf  P}^\vee$ is identified with 
${\mathbf P}(T^*{\mathbf P}^\vee)$.

The conormal bundle
$L_Q=T^*_Q({\mathbf P}\times {\mathbf P}^\vee)
\subset (T^*{\mathbf P}\times T^*{\mathbf P}^\vee)|_Q$
is a line bundle.
Since $1\in {\rm End}(V)
=V^\vee\otimes V$ regarded as
a global section of
${\mathcal O}(1)\boxtimes{\mathcal O}(1)$
is the bilinear form defining $Q
\subset {\mathbf P}\times {\mathbf P}^\vee$,
the morphism
$N_{Q/({\mathbf P}\times {\mathbf P}^\vee)}
\to \Omega^1_{({\mathbf P}\times {\mathbf P}^\vee)/
{\mathbf P}^\vee}
\otimes_{{\mathcal O}_{{\mathbf P}\times {\mathbf P}^\vee}}
{\mathcal O}_Q
=
\Omega^1_{{\mathbf P}/k}
\otimes_{{\mathcal O}_{\mathbf P}}
{\mathcal O}_Q$
defines a tautological sub invertible sheaf
on $Q={\mathbf P}(T^*{\mathbf P})$.
In other words,
the tautological sub line bundle
$L
\subset T^*{\mathbf P}\times_{\mathbf P}Q$ 
is the image of $L_Q$ by
the first projection
${\rm pr}_1\colon
(T^*{\mathbf P}\times T^*{\mathbf P}^\vee)|_Q
\to 
T^*{\mathbf P}\times_{\mathbf P}Q$.
By symmetry,
the image by
the second projection
equals
the tautological sub line bundle $L^\vee$
on $Q={\mathbf P}(T^*{\mathbf P}^\vee)$.

Since the conormal bundle $L_Q$
is the kernel of the surjection
$(T^*{\mathbf P}\times
T^*{\mathbf P}^\vee)|_Q\to T^*Q$,
the intersection
$q^\circ T^*{\mathbf P}
\cap q^{\vee \circ}T^*{\mathbf P}^\vee
\subset T^*Q$
equals the image of
the tautological bundle
$L\subset
T^*{\mathbf P}\times_{\mathbf P}Q$.
By symmetry, the intersection also 
equals the image of
the tautological bundle
$L^\vee\subset
T^*{\mathbf P}^\vee\times_{{\mathbf P}^\vee}Q$.

Let $C\subset T^*{\mathbf P}$ denote
a closed conical subset.
We define the Legendre transform
$C^\vee\subset T^*{\mathbf P}^\vee$
to be $q^\vee_\circ q^\circ C$.
We consider projectivizations
${\mathbf P}(C)
\subset {\mathbf P}(T^*{\mathbf P})$
and
${\mathbf P}(C^\vee)
\subset {\mathbf P}(T^*{\mathbf P}^\vee)$
as closed subsets of
$Q$.
Let $C^+\subset T^*{\mathbf P}$ denote
the union of $C$ and the $0$-section.

\begin{pr}\label{prL}
Let $C\subset T^*{\mathbf P}$
be a closed conical subset.
Let 
$E={\mathbf P}(C)
\subset 
Q={\mathbf P}(T^*{\mathbf P})$
be the projectivization.
Let $L_Q=T^*_Q({\mathbf P}\times {\mathbf P}^\vee)
\subset
(T^*{\mathbf P}\times T^*{\mathbf P}^\vee)|_Q$
be the conormal line bundle.

{\rm 1.}
The projectivization
$E={\mathbf P}(C)
\subset Q$
is the complement of the largest open subset
where $(q,q^\vee)$ is $C$-acyclic.


{\rm 2.}
The Legendre transform
$C^\vee$
equals the union of the image of
$L|_E
\subset q^\circ T^*{\mathbf P}
\cap q^{\vee\circ}T^*{\mathbf P}^\vee
\subset T^*Q$
and its base.
%
We have ${\mathbf P}(C)
={\mathbf P}(C^\vee)$.

{\rm 3.}
We have $C^{\vee\vee}
\subset C^+$.
\end{pr}

\proof{
1.
The kernel 
${\rm Ker}((T^*{\mathbf P}\times
T^*{\mathbf P}^\vee)|_Q\to T^*Q)$
equals the conormal bundle
$L_Q$ and the first projection
induces an isomorphism
$L_Q\to L\subset T^*{\mathbf P}\times_{\mathbf P}Q$
to the tautological bundle.
By this isomorphism,
the intersection
$(C\times T^*{\mathbf P}^\vee)|_Q
\cap L_Q$ is identified with
$q^*C\cap L$.
Hence $(q,q^\vee)$
is $C$-acyclic on $U\subset Q$
if and only if the restriction
$(q^*C\cap L)|_U$ is a subset of 
the $0$-section.
Since the projectivization
$E={\mathbf P}(C)\subset
{\mathbf P}(T^*{\mathbf P})$
equals
${\mathbf P}(q^* C\cap L)\subset
{\mathbf P}(L)=Q$,
the assertion follows.

2.
The Legendre transform $C^\vee$
is the image of 
the intersection
$q^\circ C\cap
q^{\vee \circ}T^*{\mathbf P}^\vee
\subset T^*Q$
by $q^{\vee \circ}T^*{\mathbf P}^\vee
\to T^*{\mathbf P}^\vee$.
We identify the intersection
$q^\circ T^*{\mathbf P}
\cap q^{\vee\circ}T^*{\mathbf P}^\vee
\subset T^*Q$
with the tautological line bundle
$L\subset T^*{\mathbf P}\times_{\mathbf P}Q$.
Then, the intersection
$q^\circ C\cap 
q^{\vee \circ}T^*{\mathbf P}^\vee
\subset T^*Q$
is identified with
$q^* C\cap L$.
Since the projectivization
$E={\mathbf P}(C)\subset
{\mathbf P}(T^*{\mathbf P})$
equals
${\mathbf P}(q^* C\cap L)\subset
{\mathbf P}(L)=Q$,
the closed conical subset $q^* C\cap L$
equals $L|_E$ 
outside the $0$-section.

%
%
Since 
$L|_E$ is identified with
$L^\vee|_E$ inside 
$T^*Q$,
we have
${\mathbf P}(C^\vee)=
{\mathbf P}(L^\vee|_E)
=E\subset
{\mathbf P}(T^*{\mathbf P}^\vee)=Q$.

3.
By 2 and symmetry,
we have
${\mathbf P}(C)
={\mathbf P}(C^{\vee})
={\mathbf P}(C^{\vee\vee})$.
Hence we have
$C^{\vee\vee}\subset C^+$.
\qed

}

\begin{lm}\label{lmLeg}
Let $C \subset T^*{\mathbf P}$
be a closed conical subset and
$E={\mathbf P}(C)
\subset Q=
{\mathbf P}(T^*{\mathbf P})$
be its projectivization.
Let 
$$\xymatrix{
{\mathbf P} &
Q\ar[l]_q\ar[d]_{q^\vee}&
Q_W\ar[l]_{h'}\ar[d]_{q^\vee_W}
\ar[rd]^{f'}&
\\
&{\mathbf P}^\vee
&W\ar[l]_{h}\ar[r]^f&
Y}$$ be a commutative
diagram of smooth schemes over $k$
with cartesian square.
Let $C^{\vee+}
=C^\vee\cup  T^*_{{\mathbf P}^\vee}
{\mathbf P}^\vee
\subset T^*{\mathbf P}^\vee$
be the union of the Legendre transform
with the $0$-section and
suppose that $(h,f)$ is $C^{\vee+}$-acyclic.

{\rm 1.}
The morphism $f$ is smooth
on 
$W$
and the pair
$(qh',f')$ is
$C$-acyclic on 
the complement $Q_W\sm E_W$.

{\rm 2.}
On a neighborhood of $E_W$,
the pair $(qh',f')$ is
$T^*{\mathbf P}$-acyclic.
\end{lm}

\proof{
1.
Since $(h,f)$ is
$C^{\vee+}$-acyclic and 
$C^{\vee+}$ contains the $0$-section,
the morphism $f\colon W\to Y$ is smooth.
By Proposition \ref{prL}.1,
$(q,q^\vee)$ is $C$-acyclic
outside $E$.
Hence
$(qh',q^\vee_W)$ is $C$-acyclic
outside $E_W$ by 
Lemma \ref{lmuniCacy}
and 
$(qh',f')$ is $C$-acyclic
outside $E_W$ by 
Lemma \ref{lmCacych}.

2.
By the description of $C^\vee$
in Proposition \ref{prL}.2
and by the open condition Lemma \ref{lmCtrans}.4
and Lemma \ref{lmCacycsm}.1,
the $C^{\vee}$-acyclicity of
$(h,f)$ implies
the $T^*{\mathbf P}$-acyclicity of
$(qh',f')$
on a neighborhood 
of $E_W\subset Q_W$.
\qed

}

\medskip

We define the naive Radon transform
$R{\mathcal F}$ to be
$q^\vee_*q^*{\mathcal F}$ and
the naive inverse Radon transform
$R^\vee{\mathcal G}$ to be
$q_*q^{\vee*}{\mathcal G}$.
Since we use only the naive Radon transform,
we drop the adjective `naive' in the sequel.
We will refine in Proposition \ref{prSSR},
after studying the difference between
$R^\vee R{\mathcal F}$
and ${\mathcal F}$,
the following elementary property.

\begin{lm}\label{lmL}
Assume that ${\mathcal F}$ is micro supported on $C$.

{\rm 1.}
The Radon transform
$R{\mathcal F}$ is micro supported on $C^\vee$.

{\rm 2.}
$R^\vee R{\mathcal F}$ is micro supported on $C^+$.
\end{lm}

\proof{
1.
By Lemma \ref{lmmsloc}.1 and 
Lemma \ref{lmCbarms},
the Radon transform
$R{\mathcal F}=q^\vee_*q^*{\mathcal F}$
is micro supported
on the Legendre transform
$C^\vee=q^\vee_\circ q^\circ C$.

2.
By 1,
$R^\vee R{\mathcal F}$ is micro supported on 
$C^{\vee\vee}\subset C^+$.
\qed
}

\begin{lm}\label{lmRn}
We consider the commutative diagram
$$
\xymatrix{
Q\times_{{\mathbf P}^\vee}Q
\ar[r]^-i\ar[rd]_{q\times q}&
{\mathbf P}\times
{\mathbf P}^\vee\times {\mathbf P}
\ar[d]^{{\rm pr}_{13}}
\\
&
{\mathbf P}\times
{\mathbf P}
&
{\mathbf P}\ar[l]_-{\delta_{\mathbf P}}
}
$$
where 
$\delta_{\mathbf P}
\colon {\mathbf P}\to {\mathbf P}\times
{\mathbf P}$
is the diagonal immersion.

{\rm 1.}
The closed immersion $i=((q,q^\vee),(q^\vee,q))$
induces isomorphisms
\begin{equation}
R^s(q\times q)_*\Lambda_{
Q\times_{{\mathbf P}^\vee}Q}
\to 
\begin{cases}
\Lambda(-t)[-2t]
&\text{if }s=2t
\text{ and }0\leqslant t\leqslant n-2,\\
\delta_{{\mathbf P}*}\Lambda(-(n-1))[-2(n-1)]
&\text{if }s=2(n-1),
\\0&\text{if otherwise.}
\end{cases}
\label{eqPQP}
\end{equation}

{\rm 2.}
Let 
$p\colon {\mathbf P}\to {\rm Spec}\, k$
and $p^\vee\colon {\mathbf P}^\vee\to {\rm Spec}\, k$
denote the projections.
Then, we have a distinguished triangle
\begin{equation}
\to 
\tau_{\leqq 2(n-2)}
(p\times p)^*
p^\vee_*\Lambda
\to
(q\times q)_*\Lambda_{
Q\times_{{\mathbf P}^\vee}Q}
\to
\delta_{{\mathbf P}*}
\Lambda_{\mathbf P}(n-1)[2(n-1)]\to.
\label{eqppL}
\end{equation}
\end{lm}

\proof{
1.
The immersion $i$ induces a morphism
\begin{equation}
(p\times p)^*
p^\vee_*\Lambda
={\rm pr}_{13*}\Lambda
\to 
(q\times q)_*\Lambda_{
Q\times_{{\mathbf P}^\vee}Q
}
\label{eqQQ}
\end{equation}
and we have isomorphisms
$R^sp^\vee_*\Lambda
\to 
\Lambda(-t)[-2t]$
for $s=2t, 0\leqslant t\leqslant n$
and 
$R^sp^\vee_*\Lambda
=0$ otherwise.
The restriction of
the closed immersion $
i\colon Q\times_{{\mathbf P}^\vee}Q
\to
{\mathbf P}\times 
{\mathbf P}^\vee
\times {\mathbf P}$
on the diagonal
${\mathbf P}\subset
{\mathbf P}\times {\mathbf P}$
is the sub ${\mathbf P}^{n-1}$-bundle $Q
\subset {\mathbf P}\times{\mathbf P}^\vee$.
On the complement
${\mathbf P}\times {\mathbf P}
\sm {\mathbf P}$,
$Q\times_{{\mathbf P}^\vee}Q$
is a sub
${\mathbf P}^{n-2}$-bundle.
Hence (\ref{eqQQ}) induces an isomorphism
of cohomology sheaves except
for degree $s=2(n-1)$
and induces an isomorphism
$\delta_{{\mathbf P}*}
p^*R^{2(n-1)}p^\vee_*\Lambda
\to 
R^{2(n-1)}(q\times q)_*\Lambda_{
Q\times_{{\mathbf P}^\vee}Q}$.

2.
This follows from the isomorphisms (\ref{eqPQP}).
\qed

}
\medskip

Next, we consider the diagram
\begin{equation}
\begin{CD}
{\mathbf P}@<{{\rm pr}_1}<<
{\mathbf P}\times {\mathbf P}
@<{q\times q}<<
Q\times_{{\mathbf P}^\vee}Q
\\
@.@V{{\rm pr}_2}VV@.\\
@.{\mathbf P}.@.
\end{CD}
\label{eqRRv}
\end{equation}

\begin{pr}\label{prRn}
{\rm 1.}
We have a canonical isomorphism
\begin{equation}
R^\vee R{\mathcal F}
\to
R{\rm pr}_{2*}\bigl({\rm pr}_1^*{\mathcal F}
\otimes R(q\times q)_*\Lambda_{
Q\times_{{\mathbf P}^\vee}Q}\bigr).
\label{eqRRF}
\end{equation}

{\rm 2.}
The isomorphism {\rm (\ref{eqRRF})}
induces a
distinguished triangle
$$
\to
\bigoplus_{q=0}^{n-2}
p^*p_*{\mathcal F}(-q)[-2q]
\to
R^\vee R{\mathcal F}
\to 
{\mathcal F}(-(n-1))[-2(n-1)]
\to.
$$
\end{pr}

\proof{
1.
By the cartesian diagram
$$\begin{CD}
{\mathbf P}@<q<<
Q@<{{\rm pr}_1}<<
Q\times_{{\mathbf P}^\vee}Q
\\
@.@V{q^\vee}VV@VV{{\rm pr}_2}V\\
@.{{\mathbf P}^\vee}@<{q^\vee}<<Q\\
@.@.@VVqV\\
@.@.{\mathbf P},
\end{CD}$$
we have 
$R^\vee R{\mathcal F}=
Rq_*q^{\vee*}Rq^\vee_*q^*{\mathcal F}$.
By the proper base change theorem,
we have a canonical isomorphism
$Rq_*q^{\vee*}Rq^\vee_*q^*{\mathcal F}
\to 
R(q\circ {\rm pr}_2)_*
(q\circ {\rm pr}_1)^*{\mathcal F}$.
In the notation of (\ref{eqRRv}),
the latter is identified with
$R({\rm pr}_2\circ (q\times q))_*
({\rm pr}_1\circ (q\times q))^*{\mathcal F}$.
This is identified with
$R{\rm pr}_{2*}({\rm pr}_{1*}{\mathcal F}
\otimes 
R(q\times q)_*\Lambda_{Q\times_{{\mathbf P}^\vee}Q})$ by the projection formula.

2.
This follows from the isomorphism (\ref{eqRRF})
and the distinguished triangle (\ref{eqppL}).
}

\begin{pr}\label{prSSR}
For a sheaf ${\mathcal F}$ on ${\mathbf P}$
and a closed conical subset
$C\subset T^*{\mathbf P}$,
we have implications
{\rm (1)}$\Rightarrow${\rm (2)}$\Rightarrow${\rm (3)}$\Rightarrow${\rm (4)}.

{\rm (1)}
${\mathcal F}$ is micro supported on $C$.

{\rm (2)}
$q^{\vee}$
is universally $q^*{\mathcal F}$-acyclic
 outside $E={\mathbf P}(C)$.

{\rm (3)}
The Radon transform
$R{\mathcal F}$ is micro supported on $C^{\vee+}$.

{\rm (4)}
${\mathcal F}$ is micro supported on $C^+$.
\end{pr}

%

\proof{
{\rm (1)}$\Rightarrow${\rm (2)}:
The pair $(q,q^\vee)$
of $q\colon Q\to {\mathbf P}$
and $q^\vee\colon Q\to {\mathbf P}^\vee$ is 
$C$-acyclic
outside $E={\mathbf P}(C)$
by Proposition \ref{prL}.1.
Hence (1) implies that
$q^\vee$ is universally
$q^*{\mathcal F}$-acyclic on the complement  
$Q\sm E$ by 
Lemma \ref{lmms}.4.

{\rm (2)}$\Rightarrow${\rm (3)}:
Assume that a pair of morphisms $h\colon W\to {\mathbf P}^\vee,
f\colon W\to Y$ is $C^{\vee+}$-acyclic
and show that
$f$ is $h^*R{\mathcal F}$-acyclic.
We consider the commutative diagram
$$\xymatrix{
{\mathbf P} &
Q\ar[l]_q\ar[d]_{q^\vee}&
Q_W\ar[l]_{h'}\ar[d]_{q^\vee_W}
\ar[rd]^{f'}&
\\
&{\mathbf P}^\vee
&W\ar[l]_{h}\ar[r]^f&
Y}$$
with cartesian square as in Lemma \ref{lmLeg}.
Since $q^\vee$ is proper,
it suffices to show that
$f'$ is $h'^*q^*{\mathcal F}$-acyclic
by 
Lemma \ref{lmFacycg}.

By (2), 
$q^\vee_W$ is $(qh')^*{\mathcal F}$-acyclic
on 
the complement $Q_W\sm E_W$.
By Lemma \ref{lmLeg}.1,
the morphism $f\colon W\to Y$ is smooth
on 
$W$.
Hence
$f'$ is $(qh')^*{\mathcal F}$-acyclic
on 
$Q_W\sm E_W$
by \cite[Corollaire 2.7]{App}.

By Lemma \ref{lmLeg}.2,
the pair $(qh',f')$ is
$T^*{\mathbf P}$-acyclic
on a neighborhood 
of $E_W\subset Q_W$.
Since ${\mathcal F}$ is micro supported on 
$T^*{\mathbf P}$
by Lemma \ref{lmms}.2,
$f'$ is $(qh')^*{\mathcal F}$-acyclic
on a neighborhood of $E_W$.
Thus 
$f'$ is $(qh')^*{\mathcal F}$-acyclic
as required.

(3)$\Rightarrow$(4)
By (3) and (1)$\Rightarrow$(3),
$R^\vee R{\mathcal F}$ is micro supported on $(C^{\vee+})
^{\vee+}=
C^+$.
Since $p^*p_*{\mathcal F}$
is micro supported
on the 0-section
$T^*_{\mathbf P}{\mathbf P}$,
by the distinguished triangle in
Proposition \ref{prRn}.2,
${\mathcal F}$  is also micro supported on $C^+$.
\qed

}

\begin{cor}\label{corSSR}
Let ${\mathcal F}$ be a 
sheaf on ${\mathbf P}$.
Let $E\subset Q={\mathbf P}(T^*{\mathbf P})$ be 
the complement of
the largest
open subset on which
$q^\vee$
is universally $q^*{\mathcal F}$-acyclic.
Then the closed conical
subset $C\subset T^*{\mathbf P}$
corresponding to the base
$B={\rm supp}\, {\mathcal F}$ 
and the projectivization $E
\subset {\mathbf P}(T^*{\mathbf P})$
is the singular support $SS{\mathcal F}$
of ${\mathcal F}$.
\end{cor}

\proof{
By Proposition \ref{prSSR} (2)$\Rightarrow$(4),
${\mathcal F}$ is micro supported on $C^+$.
Hence
${\mathcal F}$ is micro supported on $
C=C^+|_B$ by 
Lemma \ref{lmCA}.1.

Assume that ${\mathcal F}$ is micro supported on $C'
\subset T^*{\mathbf P}$.
Then, by Proposition \ref{prSSR} (1)$\Rightarrow$(2),
we have ${\mathbf P}(C')\supset E={\mathbf P}(C)$
since $E$ is the smallest.
Since the base of $C'$ contains
$B={\rm supp}\, {\mathcal F}$
as a subset,
we have $C'\supset C$.
\qed

}

\begin{thm}\label{thmSS}
Let $X$ be a smooth scheme
over a field $k$
and
${\mathcal F}$ be a sheaf.
Then the singular support
$SS{\mathcal F}$ of ${\mathcal F}$ exists.
\end{thm}

\proof{
By Lemma \ref{lmSSloc},
the assertion is local on $X$.
Hence, we may assume that $X$ is affine
and is a closed subscheme of
${\mathbf A}^n$.
By Proposition 
\ref{prSSXP},
we may assume that
$X={\mathbf A}^n$.
Further by Lemma \ref{lmSSloc},
we may assume that
$X={\mathbf P}^n$.
This case is proved in Corollary \ref{corSSR}.
\qed

}

\begin{cor}\label{corbc}
Let $X$ be a smooth scheme
over a field $k$
and
${\mathcal F}$ be a sheaf.
Let $k'$ be an extension of $k$
and ${\mathcal F}'$ be the pull-back of
${\mathcal F}$ on the base change $X'=X_{k'}$.
Then, we have
$SS {\mathcal F}'\subset 
(SS {\mathcal F})_{k'}$
and if $k'$ is an algebraic extension of $k$,
we have
$SS {\mathcal F}'=
(SS {\mathcal F})_{k'}$.
\end{cor}

\proof{
We may assume that $X={\mathbf P}$ is a projective
space
and it suffices to show
the inequality
$E'\subset E_{k'}$
or the equality
$E'=E_{k'}$
of the projectivizations
$E={\mathbf P}(SS {\mathcal F})$
and 
$E'={\mathbf P}(SS {\mathcal F}')$.
Since the complement
$Q\sm E$ is the largest open subset
on which $p^\vee\colon Q\to {\mathbf P}^\vee$
is universally $p^*{\mathcal F}$-acyclic
and similarly for $Q_{k'}\sm E'$,
we have $(Q\sm E)_{k'}\subset Q_{k'}\sm E'$.
Further if $k'$ is an algebraic extension,
we have $(Q\sm E)_{k'}= Q_{k'}\sm E'$.
\qed

}

\begin{thm}{\rm (\cite[1.3]{SS})}
\label{thmSSc}
Let $X$ be a smooth scheme
over a field $k$.
If ${\mathcal F}$ is a constructible sheaf,
then every irreducible component
of $SS{\mathcal F}$ has the same dimension
as $X$.
\end{thm}

\setcounter{section}{99}

\setcounter{section}5

\section{Holonomic sheaves
are constructible}

\begin{df}\label{dfhol}
Let $X$ be a smooth scheme
over a field $k$ and
let ${\mathcal F}$ be a sheaf on $X$.
We say that ${\mathcal F}$ is {\em holonomic}
if the following conditions 
{\rm (1)} and {\rm (2)}
are satisfied:

{\rm (1)} 
There exists a closed conical
subset $C\subset T^*X$
such that $\dim C\leqq \dim X$
and ${\mathcal F}$ is
micro supported on $C$.

{\rm (2)}
For every geometric point
$x$,
the stalk ${\mathcal F}_x$
is constructible.
\end{df}

Since the singular support $SS{\mathcal F}$ exists
in general,
condition (1) is equivalent to
$\dim SS{\mathcal F}\leqq \dim X$.
To emphasize that we also impose condition (2),
we say that ${\mathcal F}$ is {\em quasi-holonomic}
if only the first condition
{\rm (1)}
is satisfied,
although we will not use this notion.

\begin{lm}\label{lmhol}
Let $X$ be a smooth scheme
over a field $k$ and
let ${\mathcal F}$ be a sheaf on $X$.

{\rm 1.}
Assume that ${\mathcal F}$ is holonomic.
Then, there exists a dense
open scheme $U\subset X$ 
such that
${\mathcal F}_U$ is locally constant.

{\rm 2.}
Let $(U_i\to X)_{i\in I}$
be a finite family of \'etale morphism
and suppose that
it is an \'etale covering of
an open neighborhood of 
the support of ${\mathcal F}$.
Then the following conditions
are equivalent:

{\rm (1)}
${\mathcal F}$ is holonomic.

{\rm (2)}
${\mathcal F}_{U_i}$ is holonomic for every $i\in I$.

{\rm 3.}
Let 
$i\colon X\to P$ be a closed immersion
of smooth schemes over $k$.
Then, the following conditions
are equivalent:

{\rm (1)}
${\mathcal F}$ is holonomic.

{\rm (2)}
$i_*{\mathcal F}$ is holonomic.

{\rm 4.}
Let ${\mathcal F}'\to {\mathcal F}\to {\mathcal F}''\to$
be a distinguished triangle of sheaves
on $X$.
If ${\mathcal F}$ and ${\mathcal F}'$ are
holonomic,
then ${\mathcal F}''$ is holonomic.
\end{lm}

\proof{
1.
There exists a dense
open scheme $U\subset X$ 
such that
${\mathcal F}_U$ is micro supported
on the $0$-section
by Lemma \ref{lmU} 
and 
an easy case of
Lemma \ref{lmmsloc}.1.
By Lemma \ref{lmms}.1
(1)$\Rightarrow$(2),
${\mathcal F}_U$ is locally constant.

2.
By Lemma \ref{lmCA}.1,
we may assume that
$(U_i\to X)_{i\in I}$ is an \'etale covering.

(1)$\Rightarrow$(2):
If ${\mathcal F}$ is micro supported
on $C\subset T^*X$
and if $\dim C\leqq \dim X$,
then
${\mathcal F}_{U_i}$ is micro supported
on $C_{U_i}\subset T^*U_i$ 
by 
Lemma \ref{lmmsloc}.1
and $\dim C_{U_i}\leqq \dim U_i$
for every $i\in I$.

(2)$\Rightarrow$(1):
Suppose 
${\mathcal F}_{U_i}$ is micro supported
on $C_i\subset T^*U_i$
and $\dim C_i\leqq \dim U_i$
for every $i\in I$.
Then,
${\mathcal F}$ is micro supported
on $C=\overline
{\bigcup_{i\in I}{\rm Im}(C_i)}
\subset T^*X$
by Lemma \ref{lmmsloc}.2.
Since $I$ is finite,
we have
$C
=
\bigcup_{i\in I}\overline
{{\rm Im}(C_i)}$
and $\dim C
=\max_{i\in I}\dim C_i\leqq \dim X$.

3.
(1)$\Rightarrow$(2):
If ${\mathcal F}$ is micro supported
on $C\subset T^*X$
and if $\dim C\leqq \dim X$,
then
$i_*{\mathcal F}$ is micro supported
on $i_\circ C\subset T^*P$
by Lemma \ref{lmCA}.3
and $\dim i_\circ C\leqq \dim P$.

(2)$\Rightarrow$(1):
Suppose that
the singular support $C_P=SS i_*{\mathcal F}$
satisfies
$\dim C_P\leqq \dim P$.
By Lemma \ref{prSSXP}.2,
the singular support
$C=SS{\mathcal F}$ exists
and $C_P=i_\circ C$.
Hence, we have
$\dim C=\dim C_P-{\rm codim}_PX
\leqq \dim X$.

4.
Suppose that 
${\mathcal F}$ and ${\mathcal F}'$
are micro supported
on $C$ and $C'$
respectively and
that $\dim C,\dim C'\leqq \dim X$.
By Lemma \ref{lmmsUi},
${\mathcal F}''$ is micro supported
on $C''=C\cup C'$
and we have
$\dim C''\leqq \dim X$.
\qed

}
\medskip

Theorem \ref{thmSSc}
states that every constructible
sheaf is holonomic as
the title of the article \cite{SS} says.
We have a partial converse.

\begin{thm}\label{thmhol}
Let $X$ be a smooth scheme
of finite type over a 
field $k$
and ${\mathcal F}$ be a sheaf
on $X$.
If ${\mathcal F}$ is holonomic,
then ${\mathcal F}$ is constructible.
\end{thm}

Without condition (2)
in Definition \ref{dfhol},
we have an obvious counterexample
where $X={\rm Spec}\, k$.
Condition (2) is not enough for the converse
since
${\mathcal F}=\bigoplus_{x: \text{\rm closed points of }X}
\Lambda_x$ satisfies (2) but
its micro support is the whole $T^*X$.

\proof{
We prove the theorem
first assuming that $k$ is perfect 
by induction
on $\dim X$.
If $X=\varnothing$,
the assertion is clear.
By Lemma \ref{lmhol}.1,
the largest open subset
$U\subset X$
on which ${\mathcal F}$ is constructible
is nonempty.
To prove the equality $X=U$
by contradiction,
we assume that the complement
$Z=X\sm U$
is not empty.
Since $k$ is assumed perfect,
after shrinking $X$ if necessary,
we may assume that $Z$
is smooth over $k$.

Let $j\colon U\to X$
be the open immersion,
$i\colon Z\to X$
be the closed immersion
and consider
the distinguished
triangle
$j_!j^*{\mathcal F}\to {\mathcal F}
\to i_*i^*{\mathcal F}\to$.
Since $j_!j^*{\mathcal F}$
is constructible,
it suffices to show 
that $i_*i^*{\mathcal F}$ is constructible.
By Theorem \ref{thmSSc},
$j_!j^*{\mathcal F}$ is holonomic.
Hence by Lemma \ref{lmhol}.4,
$i_*i^*{\mathcal F}$ is also holonomic.
Since replacing ${\mathcal F}$
by $i_*i^*{\mathcal F}$ does not
change $U$,
we may assume that
${\mathcal F}=i_*{\mathcal F}_Z$
for ${\mathcal F}_Z=i^*{\mathcal F}$.
Since $i_*{\mathcal F}_Z$ is holonomic,
${\mathcal F}_Z$ is also holonomic 
by Lemma \ref{lmhol}.3.
Since $\dim Z<\dim X$,
${\mathcal F}_Z$ is constructible.
This contradict
the assumption that $U$ is the largest.

We show the general case.
Let $k'$ be a perfect closure of $k$.
Then, by Corollary \ref{corbc},
the pull-back ${\mathcal F}'$
on the base change $X_{k'}$ is holonomic.
Hence by the perfect case,
${\mathcal F}'$ is constructible.
Since the canonical morphism
$X_{k'}\to X$ is a homeomorphism,
${\mathcal F}$ is also constructible.
\qed

}
\medskip

We may extend Definition \ref{dfhol}
to sheaves on singular schemes as follows.

\begin{lm}\label{lmhols}
Let $X$ be a scheme of finite
type over a field $k$ and
${\mathcal F}$ be a sheaf on $X$.

{\rm 1.}
Let $i\colon X\to P$
and $i'\colon X\to Q$
be closed immersions to
smooth schemes
over $k$.
Then, the following conditions
are equivalent:

{\rm (1)}
$i_*{\mathcal F}$ is holonomic.

{\rm (2)}
$i'_*{\mathcal F}$ is holonomic.

{\rm 2.}
The following conditions
are equivalent:

{\rm (1)}
For every \'etale morphism
$U\to X$
and every closed immersion
$i\colon U\to P$ to a smooth scheme
over $k$,
$i_*({\mathcal F}|_U)$ is holonomic.

{\rm (2)}
There exists an \'etale covering
$(U_j\to X)_{j\in J}$
and a family of closed immersions
$i'_j\colon U_j\to Q_j$ to smooth schemes
over $k$,
such that $i'_{j*}({\mathcal F}|_{U_j})$ is holonomic
for every $j$.
\end{lm}

\proof{
1.
It suffices to show
(1)$\Rightarrow$(2).
By replacing $Q$ by $P\times Q$,
we may assume that there exists
a smooth morphism $Q\to P$
compatible with $i$ and $i'$.
Since
$i'$ induces a section
of a smooth morphism
$Q\times_PX\to X$
and since the question is local on $X$ by
Lemma \ref{lmhol}.2,
we may assume that there exists
a closed subscheme $Q'$ of $Q$
containing $X$ and \'etale over $P$
by \cite[Chapitre 0$_{\rm IV}$ (15.1.16)
a)$\Leftrightarrow$c)]{EGA4}
and \cite[Th\'eor\`eme (17.6.1)
a)$\Leftrightarrow$c$'$)]{EGA4}.
By Lemma \ref{lmhol}.3,
replacing $Q$ by $Q'$,
we may assume that there exists
an \'etale morphism $h\colon Q\to P$
compatible with $i$ and $i'$.
Since $X\to Q$ induces an open immersion
$X\to X\times_PQ$,
after shrinking $Q$, we may assume
that $X\to X\times_PQ$ is an isomorphism.
If $i_*{\mathcal F}$ is micro supported on
$C\subset T^*P$ and $\dim C\leqq \dim P$,
then
$i'_*{\mathcal F}$ is micro supported on
$h^*C\subset T^*Q$ 
by Lemma \ref{lmmsloc}.1
and $\dim h^*C
=\dim C\leqq \dim P=\dim Q$.

2.
The implication (1)$\Rightarrow$(2) is clear.

(2)$\Rightarrow$(1): 
Since $X$ is quasi-compact,
we may assume that $J$ is finite.
By \cite[Corollaire (18.4.7)]{EGA4},
for each $j\in J$,
there exist an open covering
$(U_{jh}\to U\times_XU_j)_{h\in H_j}$,
closed immersions
$i_{jh}\colon U_{jh}\to P_{jh}$
to schemes \'etale over
$P$ and
closed immersions
$i'_{jh}\colon U_{jh}\to Q_{jh}$
to schemes \'etale over
$Q_j$.
Then since $U_{jh}\to U_j\times_{Q_j}Q_{jh}$
is an open and closed immersion,
$i'_{jh*}({\mathcal F}|_{U_{jh}})$ is holonomic
for every $j,h$
by Lemma \ref{lmhol}.2.
Hence
$i_{jh*}({\mathcal F}|_{U_{jh}})$
is also holonomic by 1.
Since $(P_{jh}\to P)_{j\in J,h\in H_j}$
is an \'etale covering
of the support of an open neighborhood
of the image $i(U)$,
$i_*({\mathcal F}|_U)$
is holonomic by
Lemma \ref{lmhol}.2.
\qed

}

\begin{df}\label{dfhols}
Let ${\mathcal F}$ be a sheaf on 
a scheme $X$ of finite type over 
a field $k$.
We say that 
${\mathcal F}$ is {\em holonomic}
if 
the equivalent conditions
in Lemma {\rm \ref{lmhols}.2}
are satisfied.
\end{df}

\begin{pr}\label{prhol}
Let ${\mathcal F}$ be a sheaf on 
a scheme $X$ of finite type over 
a field $k$.
Then, the following conditions are equivalent:

{\rm (1)}
${\mathcal F}$ is constructible.

{\rm (2)}
${\mathcal F}$ is holonomic.
\end{pr}

\proof{
Since the question is local on $X$,
we may assume that
there exists a closed immersion
$i\colon X\to P$
to a smooth scheme over $k$.
Then $i_*{\mathcal F}$ is constructible
if and only if it is holonomic
by Theorem \ref{thmSSc}
and Theorem \ref{thmhol}.
\qed

}
\medskip

We have the following corollary
since the corresponding
properties are known for
constructible sheaves
\cite{Fu}.

\begin{cor}\label{corhol}
Let ${\mathcal F}$ be a sheaf on 
a scheme $X$ of finite type over 
a field $k$.

{\rm 1.}
Let $f\colon X\to Y$ be
a morphism of schemes 
of finite type over
$k$.
If ${\mathcal F}$ is holonomic,
then
$f_*{\mathcal F}$ is holonomic.
Further if $f$ is separated,
then
$f_!{\mathcal F}$ is holonomic.

{\rm 2.}
Let $h\colon W\to X$ be
a morphism of schemes 
of finite type over
$k$.
If ${\mathcal F}$ is holonomic,
then
$h^*{\mathcal F}$ is holonomic.
Further if $h$ is separated,
then
$h^!{\mathcal F}$ is holonomic.
\end{cor}


\begin{thebibliography}{99}

\bibitem{Artin}
M.~Artin,
{\em Faisceaux constructibles, cohomologie d'une courbe alg\'ebrique},
SGA 4 Expos\'e IX. Th\'eorie des Topos 
et Cohomologie \'Etale des Sch\'emas, 
Lecture Notes in Mathematics, vol. 305, 
pp. 1--42 (1973)

\bibitem{AL}
-----,
{\em Morphismes acycliques},
SGA 4 Expos\'e XV. Th\'eorie des Topos 
et Cohomologie \'Etale des Sch\'emas, 
Lecture Notes in Mathematics, vol. 305, 
pp. 168--205 (1973)


\bibitem{SS}
A.~Beilinson, {\em Constructible sheaves are holonomic},
Selecta Math., 22 (4), 1797--1819 (2016).

\bibitem{TF}
P.~Deligne, {\em Th\'eor\`emes de finitude en cohomologie $\ell$-adique}, Cohomologie \'etale SGA 4$\frac12$.
Springer Lecture Notes in Math. 569, 
233--251 (1977)

\bibitem{Fu}
L.~Fu,
{\sc Etale cohomology theory},
Nankai tracts in Mathematics Vol.~13,
World Scientific, 2011,
revised ed.~2015.



\bibitem{EGA4}
A.~Grothendieck,
{\sc \'El\'ements de g\'eom\'etrie alg\'ebrique IV-4}, 
\'Etude locale des sch\'emas et
des morphismes de sch\'emas, 
Publ.\ Math.\ IHES 32 (1967).


\bibitem{HTT}
R.~Hotta, K.~Takeuchi, T.~Tanisaki,
{\sc $D$-Modules, Perverse Sheaves, and Representation Theory},
Progress in Mathematics,
Volume 236,
Birkha\"user
Boston, Basel, Berlin,
2008, translated from Japanese original
1995.


\bibitem{App}
L.~Illusie,
{\em Appendice \`a Th\'eor\`eme de finitude en cohomologie $\ell$-adique},
Cohomologie \'etale SGA 4$\frac12$. Springer Lecture Notes in Mathematics, 
vol.~569, pp.~252--261 (1977)


\bibitem{Kas}
M.~Kashiwara,
{\em On the Maximally Overdetermined System
of Linear Differential Equations, I},
Publications RIMS, Kyoto Univ.~10 (1975), 563--579

\bibitem{SKK}
M.~Sato, T.~Kawai, M.~Kashiwara,
{\it Mircofunctions and pseudo-differential
equations}, in Hyperfunctions and pseudo-differential
equations (Katata, 1971),
Lecture Notes in Math., vol 287,
Springer-Verlag, 1973,
pp. 265--529.


\end{thebibliography}
\end{document}